\documentclass[12pt]{article}
\usepackage{amsfonts}
\usepackage{latexsym,amsmath,amssymb}
\usepackage{graphics,epstopdf}
\usepackage[pdftex]{graphicx}

\newtheorem{definition}{Definition}[section]
\newtheorem{theorem}[definition]{Theorem}

\newtheorem{lemma}[definition]{Lemma}
\numberwithin{equation}{section}
\begin{document}

\begin{center}
{\Large \textbf{\ On approximation by Stancu type q-Bernstein-Schurer-Kantorovich operators }}

\bigskip

\textbf{M. Mursaleen} and \textbf{Taqseer Khan} 

Department of\ Mathematics, Aligarh Muslim University, Aligarh--202002, India%
\\[0pt]

mursaleenm@gmail.com; taqi.khan91@gmail.com;

\bigskip

\bigskip

\textbf{Abstract}
\end{center}

\parindent=8mm {\footnotesize {In this paper we introduce the Stancu type generalization of the q-Bernstein-Schurer-Kantorovich operators
 and examine their approximation properties. We investigate the convergence of our operators with the help of the Korovkin's approximation theorem
 and examine the convergence of these operators in the Lipschitz class of functions. We also investigate the approximation process for these operators through
 the statistical Korovkin's approximation theorem. Also, we present some direct theorems for these operators. Finally we introduce the bivariate
 analogue of these operators and study some results for the bivariate case.}}
\bigskip

{\footnotesize \emph{Keywords}: Stancu type q-Bernstein-Schurer-Kantorovich operators; modulus of continuity; positive
linear operators; Korovkin type approximation theorem;statistical approximation, Lipschitz class of functions.


\section{Introduction and preliminaries}
The $q$-calculus has played an inportant role in the field of approximation theorey since last three decades. In the year of 1987, A. Lupas
was the first to apply the $q$-calculus to approximation theory. He introduced the $q$-analogue of the well known Bernstein polynomials $[16]$.
Another remarkable application of the $q$-calculus advented in the year of 1997 by Phillips $[14]$. He used the
$q$-calculus to define another interesting $q$-analogue of the classical Bernstein polynomials. Ostrovska [15] obtained more results on the $q$-
Bernstein polynomials. In the sequel many researchers have studied the $q$- analogues of many well known operators like Baskakov operators,
Meyer-Konig-Zeller operators, Szas-Mirakyan operators, Bleiman, Butzer and Hahn operators (written succinctly as BBH).
Also the $q$-analogues of some integral operators like Kantorovich and Durrmeyer type were introduced and their approximation properties were studied.\\
In [20] Muraru defined the $q$- Bernstein-Schurer operators in 2011. She used the modulus of continuity to obtain the rate of convergence of the
$q$-Bernstein-Schurer operators. Recently, the Schurer modifications of some positive linear operators have been studied in [1, 2, 3].\\

Kantorovich introduced the following integral type generalization of the classical Bernstein operatos
\begin{equation*}
L_n(f;x) = (n+1)\sum_{k=0}^{n}p_{n,k}(x)\int_{\frac{k}{n+1}}^\frac{k+1}{n+1} f(t)dt,
\end{equation*}
where
\begin{equation*}
p_{n,k}(x)= \binom{n}{k}x^k (1-x)^{n-k}, ~x\in [0, 1],
\end{equation*}
is the Bernstein basis function. The approximation properties of these operators through the Korovkin's
approximation theorem were studied in [5]. \\
Dalmanoglu presented another Kantorovich type generalization of the $q$-Bernstein polynomials and studied some approximation results in [9].
Below we give some rudiments of the $q$-calculus. In [10], Radu investigated the statistical convergence properties of the Bernstein-Kantorovich
polynomials based on $q$-integers. Recently, many researchers have studied various $q$-extensions of the Kantorovich operatos in [11, 12, 13]. \\
 For any fixed real number $q>0$  and $k\in N\cup\{0\}$, the $q$-integer of $k$, denoted by $[k]_q$, is defined by
 \[ [k]_{q}=\left\{
\begin{array}{ll}
      \frac{1-q^k}{1-q},~~ & q\neq1 \\
       ~k, & q=1 \\
\end{array}
\right. \]
and the $q$-factorial $[k]_q!$ is defined as
\[ [k]_{q}!=\left\{
\begin{array}{ll}
      [k]_q[k-1]_q ...[1]_q,~~ & k\geq1 \\
       ~1, & k=0. \\
\end{array}
\right. \]
The $q$- concept can be extended to any real number $k$.
For integers $n$ and $k$ such that $0\leq k\leq n$, the $q$-analogue of the binomial coefficient is defined by
\begin{equation*}
\binom{n}{k}_q = \frac{[n]_q!}{[k]_q! [n-k]_q!}.
\end{equation*}
For the $q$- binomial coefficient the following relations hold:
\begin{eqnarray*}
  \binom{n}{k}_q&=& \binom{n-1}{k-1}_q + q^k \binom{n-1}{k}_q\\
  \binom{n}{k}_q &=& q^{n-k} \binom{n-1}{k-1}_q + \binom{n-1}{k}_q.
\end{eqnarray*}
For $x\in [0,1]$ and $m\in N^0$, the $q$-analogue of $(1+x)^n$, denoted by $(1+x)^n_q$, is defined by the polynomial
\[ (1+x)_q^n=\left\{
\begin{array}{ll}
      (1+x)(1+qx) ...(1+q^{n-1}x),~~ & n = 1, 2, 3, ... \\
       ~1, & n=0. \\
\end{array}
\right. \]
For $0<q<1, ~a>0$, the $q$-definite integral of a real valued function $f$ is defined by
\begin{equation*}
\int_{0}^{a}f(x)d_qx = (1-q)a \sum_{n=0}^{\infty}f(aq^n)q^n,~a\in R,
\end{equation*}
and over the interval $[a, b]$, $0<a <b$, it is defined by
\begin{equation*}
\int_{a}^{b}f(x)d_qx =\int_{0}^{b}f(x)d_qx -\int_{0}^{a}f(x)d_qx.
\end{equation*}
For details of these integrals one is referred to [1], [2]. In some situations the above integrals are not appropriate to obtain the $q$-analogues of
well known integrals. So we use another more general integrals, the Reimann type $q$-integrals, defined as follows:
\begin{equation*}
\int_{a}^{b}f(x)d_q^Rx = (1-q)(b-a) \sum_{s=0}^{\infty}f(a+(b-a)q^s)q^s,
\end{equation*}
where $a, b$ are such that $0\leq a<b$ and $q$ is as above. The later integrals were introduced by Gauchman [3] and Marinkovi$\acute{c}$ et al.[4].\\
In this paper, let $I$ denote the interval $[0, 1+l]$ equipped with the norm $\Vert.\Vert_{C[0, 1+l]}$, where $l\in N^0 = N\cup \{0\}$.

\section{Construction of Operators}
 In 2015, P.N. Agarwal et al. [13] introduced the following Kantorovich type generalization of the $q$-Bernstein-Schurer operators
\begin{equation}
K_{n,p}(f;q,x)=[n+1]_q\sum\limits_{k=0}^{n+p}b_{n+p,k}^q(x)q^{-k}\int_{\frac{[k]_q}{[n+1]_q}}^{\frac{[k+1]_q}{[n+1]_q}}
f(t)d_q^Rt, ~~x\in[0,1],
\end{equation}
where $b_{n+p,k}^q(x)= \binom{n+p}{k}_q x^k(1-x)_q^{n+p-k}$ is the q-Bernstein basis function. They have investigated the approximation properties of these operators using
the Korovkin's approximation theorem.
 Inspired by their work, we introduce the Stancu type
generalisation of the Bernstein-Schurer-Kantorovich operators based on q-integers as follows:
\begin{equation}\label{e2.1}
L_{n,l}^{\alpha,\beta}(f;q;x)=([n+1]_q+\beta)\sum\limits_{k=0}^{n+l}b_{n,l}^k(q;x)q^{-k}\int_{\frac{[k]_q+\alpha}{[n+1]_q+\beta}}^{\frac{[k+1]_q+\alpha}{[n+1]_q+\beta}}
f(t)d_q^Rt, ~~x\in[0,1],
\end{equation}
where $b_{n,l}^k(q;x)= \binom{n+l}{k}_q x^k(1-x)_q^{n+l-k}$ is the q-Bernstein basis function and $\alpha, \beta$ are such that ~$0<\beta\leq\alpha$.

$\alpha=0,~\beta=0$ reduce the operators (2.2) to the operators (2.1). So the newly constructed operators are a generalization of the operators in (2.1).
We shall investigate some approximation results for the operators in (2.2).
To examine the approximation results, we need the following lemmas.

\begin{lemma} Let $L_{n,l}^{\alpha,\beta}(f;q;x)$ be given by (2.1). Then the followings hold:
\begin{enumerate}
\item [(i)]$L_{n,l}^{\alpha,\beta}(1;q;x)=1$,
\item [(ii)]$L_{n,l}^{\alpha,\beta}(t;q;x)=\frac{\alpha}{[n+l]_q}+\frac{1}{([n+1]_q+\beta)[2]_q}+\frac{2q[n+l]_q}{([n+1]_q+\beta)[2]_q}x$,
\item [(iii)]$L_{n,l}^{\alpha,\beta}(t^2;q;x)=\frac{1}{([n+1]_q+\beta)^2[3]_q}+\frac{2\alpha}{([n+1]_q+\beta)^2[2]_q}+\frac{\alpha^2}{([n+1]_q+\beta)^2}
    +\frac{q[n+l]_q((3+4\alpha)+(5+4\alpha)q+4(1+\alpha)q^2)}{([n+1]_q+\beta)^2[2]_q[3]_q}x\\
    ~~~~~~~~~~~~~~~~~~+\frac{q^2[n+l]_q[n+l-1]_q(1+q+4q^2)}{([n+1]_q+\beta)^2[2]_q[3]_q}x^2$.
\end{enumerate}

\end{lemma}
\parindent=0mm\textbf{Proof}.
Before proving the above lemma, we shall first prove the following:
\begin{equation}
\sum_{k=0}^{n+l}b_{n,l}^k(q;x)q^k = 1-(1-q)[n+l]_q x,
\end{equation}
and
\begin{equation}
\sum_{k=0}^{n+l}b_{n,l}^k(q;x)q^{2k} = 1-(1-q^2)[n+l]_q x +q(1-q)^2[n+l][n+l-1]_q x^2
\end{equation}
where $b_{n,l}^k(q;x) = \binom{n+l}{k}_q x^k (1-x)_q^{n+l-k}$.
\newpage
In fact we have
\begin{eqnarray*}
\sum_{k=0}^{n+l}b_{n,l}^k(q;x)q^k &=& \sum_{k=0}^{n+l}b_{n,l}^k(q;x)\big( 1-1+q) [k]_q\big) \\
&=&\sum_{k=0}^{n+l}b_{n,l}^k(q;x)-(1-q)[n+l]_q \sum_{k=0}^{n+l}b_{n,l}^k(q;x) \frac{[k]_q}{[n+l]_q} \\
&=& 1-(1-q)[n+l]_q \sum_{k=0}^{n+l-1}\binom{n+l-1}{k}_q x^{k+1}(1-x)^{n+l-k-1}_q \\
&=& 1-(1-q)[n+l]_q x,
\end{eqnarray*}
and
\begin{eqnarray*}
\sum_{k=0}^{n+l}b_{n,l}^k(q;x)q^{2k} &=& \sum_{k=0}^{n+l}b_{n,l}^k(q;x)\bigg( (1-1+q^2)[k]_q+ q(1-q)^2[k-1]_q[k]_q\bigg) \\
&=& \sum_{k=0}^{n+l}b_{n,l}^k(q;x)-(1-q^2)[n+l]_q \sum_{k=0}^{n+l}b_{n,l}^k(q;x)\frac{[k]_q}{[n+l]_q}+ q(1-q)^2[n+l]_q [n+l-1]_q \\
&&\times\sum_{k=0}^{n+l}b_{n,l}^k(q;x)\frac{[k]_q}{[n+l]_q}\frac{[k-1]_q}{[n-1]_q} \\
&=& 1-(1-q^2)[n+l]_q \sum_{k=0}^{n+l}\binom{n+l-1}{k}_q x^k x (1-x)_q^{n+l-k-1}+q(1-q)^2 [n+l]_q[n+l-1]_q \\
&&\times\sum_{k=0}^{n+l-2}
\binom{n+l-2}{k}_q x^k x^2 (1-x)_q^{n+l-k-2}.
\end{eqnarray*}
(i)
\begin{eqnarray*}
\L_{n,l}^{\alpha,\beta}(1;q;x)&=&([n+1]_q+\beta)\sum\limits_{k=0}^{n+l}b_{n,l}^k(q;x)q^{-k}\int_{\frac{[k]_q+\alpha}{[n+1]_q+\beta}}^{\frac{[k+1]_q
+\alpha}{[n+1]_q+\beta}}1 d_q^Rt\\
&=&([n+1]_q+\beta)\sum\limits_{k=0}^{n+l}b_{n,l}^k(q;x)q^{-k} (1-q)\frac{([k+1]_q-[k]_q)}{([n+1]_q+\beta)}\sum\limits_{s=0}^{\infty}q^s \\
&=&\sum\limits_{k=0}^{n+l}b_{n,l}^k(q;x)=1.
\end{eqnarray*}
(ii)
\begin{eqnarray*}
\L_{n,l}^{\alpha,\beta}(t;q;x)&=&([n+1]_q+\beta)\sum\limits_{k=0}^{n+l}b_{n,l}^k(q;x)q^{-k}\int_{\frac{[k]_q+\alpha}{[n+1]_q+\beta}}^{\frac{[k+1]_q
+\alpha}{[n+1]_q+\beta}} td_q^Rt \\
&=&([n+1]_q+\beta)\sum\limits_{k=0}^{n+l}b_{n,l}^k(q;x)q^{-k}(1-q)\frac{([k+1]_q-[k]_q)}{([n+1]_q+\beta)}\sum\limits_{s=0}^{\infty}f
\left(\frac{[k]_q+\alpha}{[n+1]_q+\beta}+ \frac{[k+1]_q-[k]_q}{[n+1]_q+\beta}q^s\right)q^s \\
&=&\sum\limits_{k=0}^{n+l}b_{n,l}^k(q;x)(1-q)\sum\limits_{s=0}^{\infty}
\left(\frac{[k]_q+\alpha}{[n+1]_q+\beta}+ \frac{q^k q^s}{[n+1]_q+\beta}\right)q^s \\
&=&\sum\limits_{k=0}^{n+l}b_{n,l}^k(q;x)\left(\frac{[k]_q+\alpha}{([n+1]_q+\beta)}+ \frac{q^k}{([n+1]_q+\beta)}\frac{1}{[2]_q}\right)
\end{eqnarray*}
\begin{eqnarray*}
{}&=& \frac{[n+l]_q}{[n+1]_q+\beta}\sum\limits_{k=0}^{n+l}b_{n,l}^k(q;x)\frac{[k]_q+\alpha}{[n+l]_q}+\frac{1}{[2]_q([n+1]_q+\beta)}\left(1-(1-q)[n+p]_qx\right)\\
&=& \frac{[n+l]_q}{[n+1]_q+\beta}\Bigg[\sum\limits_{k=0}^{n+l}b_{n,l}^k(q;x)\frac{[k]_q}{[n+l]_q}+\sum\limits_{k=0}^{n+l}b_{n,p}^k(q;x)\frac{\alpha}{[n+l]_q}\Bigg]
+\frac{1-(1-q)[n+l]_qx}{[2]_q([n+1]_q+\beta)} \\
&=& \frac{[n+l]_q}{[n+1]_q+\beta}\Bigg[\sum\limits_{k=0}^{n+l}\binom{n+l}{k}_q x^k(1-x)^{n+l-k} \frac{[k]_q}{[n+l]_q}+ \frac{\alpha}{[n+l]_q}\Bigg]
+\frac{1-(1-q)[n+l]_qx}{[2]_q([n+1]_q+\beta)} \\
&=& \frac{[n+l]_q}{[n+1]_q+\beta}\Bigg[\sum\limits_{k=1}^{n+l-1}\binom{n+l-1}{k-1}_q x^k(1-x)^{n+l-k-1} + \frac{\alpha}{[n+l]_q}\Bigg]
+\frac{1-(1-q)[n+l]_qx}{[2]_q([n+1]_q+\beta)} \\
&=& \frac{[n+l]_q}{[n+1]_q+\beta}x + \frac{\alpha}{[n+l]_q}+\frac{1}{[2]_q([n+1]_q+\beta)}-\frac{(1-q)[n+l]_q x}{[2]_q([n+1]_q+\beta)} \\
&=& \Bigg[\frac{[n+l]_q}{[n+1]_q+\beta} -\frac{(1-q)[n+l]_q}{[2]_q([n+1]_q+\beta)}\Bigg]x+\frac{\alpha}{[n+l]_q}+\frac{1}{[2]_q([n+1]_q+\beta)} \\
&=& \frac{\alpha}{[n+l]_q}+\frac{1}{[2]_q([n+1]_q+\beta)}+\frac{2q[n+l]_q}{[2]_q([n+1]_q+\beta)}x.
\end{eqnarray*}
(iii)
\begin{eqnarray*}
\L_{n,l}^{\alpha,\beta}(t^2;q;x)&=&([n+1]_q+\beta)\sum\limits_{k=0}^{n+l}b_{n,l}^k(q;x)q^{-k}\frac{[k+1]_q-[k]_q}{([n+1]_q+\beta)}
\sum\limits_{s=0}^{\infty}\Bigg(\frac{[k]_q+\alpha}{[n+1]_q+\beta}+\frac{[k+1]_q-[k]_q}{[n+1]_q+\beta}q^s\Bigg)^2q^s \\
&=&\sum\limits_{k=0}^{n+l}b_{n,l}^k(q;x)(1-q)\sum\limits_{s=0}^{\infty}\Bigg(\frac{([k]_q+\alpha)^2}{([n+1]_q+\beta)^2}
+\frac{q^{2k}q^{2s}}{([n+1]_q+\beta)^2}+\frac{2q^k q^s([k]_q+\alpha)}{([n+1]_q+\beta)^2}\Bigg)q^s \\
&=&\sum\limits_{k=0}^{n+l}b_{n,l}^k(q;x)\Bigg(\frac{([k]_q+\alpha)^2}{([n+1]_q+\beta)^2}
+\frac{2q^k([k]_q+\alpha)}{([n+1]_q+\beta)^2(1+q)}+\frac{q^{2k}}{([n+1]_q+\beta)^2(1+q+q^2)}\Bigg) \\
&=&\sum\limits_{k=0}^{n+l}b_{n,l}^k(q;x)\frac{[k]_q^2}{([n+1]_q+\beta)^2}
+\sum\limits_{k=0}^{n+l}b_{n,l}^k(q;x)\frac{2\alpha [k]_q}{([n+1]_q+\beta)^2}+ \sum\limits_{k=0}^{n+l}b_{n,p}^k(q;x)\frac{\alpha^2}{([n+1]_q+\beta)^2} \\
&&+\sum\limits_{k=0}^{n+l}b_{n,l}^k(q;x)\frac{2q^k[k]_q}{([n+1]_q+\beta)^2 [2]_q}
+\sum\limits_{k=0}^{n+l}b_{n,l}^k(q;x)\frac{2\alpha q^k}{([n+1]_q+\beta)^2 [2]_q} \\
&&+\sum\limits_{k=0}^{n+l}b_{n,l}^k(q;x)\frac{q^{2k}}{([n+1]_q+\beta)^2[3]_q} \\
&=& \frac{[n+l]_q}{[n+1]_q+\beta}x \Bigg(\frac{1}{[n+1]_q+\beta} + q\frac{[n+l-1]_q}{[n+1]_q+\beta}x\Bigg)++\frac{2\alpha [n+l]_q}{([n+1]_q+\beta)^2}x
+\frac{\alpha^2}{([n+1]_q+\beta)^2} \\
&&+ \frac{2q [n+l]_q}{[2]_q ([n+1]_q+\beta)^2}x- \frac{2q(1-q)[n+l]_q [n+l-1]_q}{[2]_q([n+1]_q+\beta)^2}x^2
+\frac{2\alpha}{[2]_q([n+1]_q+\beta)^2} \\
&&\times\bigg(1-(1-q)[n+l]_q x\bigg)+\frac{1}{[3]_q ([n+1]_q+\beta)^2} \\
&&\times\Bigg( 1-(1-q^2)[n+l]_q x+q(1-q)^2[n+l]_q[n+l-1]_q x^2\Bigg)
\end{eqnarray*}
\begin{eqnarray*}
&=& \frac{1}{[3]_q([n+1]_q+\beta)^2} +\frac{2\alpha}{[2]_q([n+1]_q+\beta)^2}+\frac{\alpha^2}{([n+1]_q+\beta)^2} +\bigg[ \frac{[n+l]_q}{([n+1]_q+\beta)^2}
+\frac{2\alpha[n+l]_q}{([n+1]_q+\beta)^2} \\
&&+\frac{2q[n+l]_q}{[2]_q([n+1]_q+\beta)^2}-\frac{2\alpha (1-q)[n+l]_q}{[2]_q([n+1]_q+\beta)^2}-\frac{(1-q^2)[n+l]_q}{[3]_q([n+1]_q+\beta)^2} \bigg]x
+\bigg[\frac{q[n+l]_q[n+l-1]_q}{([n+1]_q+\beta)^2} \\
&&-\frac{2q(1-q)[n+l]_q[n+l-1]_q }{[2]_q([n+1]_q+\beta)^2}+\frac{q(1-q)^2[n+l]_q[n+l-1]_q}{[3]_q([n+1]_q+\beta)^2}\bigg]x^2 \\
&=& \frac{1}{[3]_q([n+1]_q+\beta)^2} +\frac{2\alpha}{[2]_q([n+1]_q+\beta)^2}+\frac{\alpha^2}{([n+1]_q+\beta)^2}+ \frac{[n+l]_q}{([n+1]_q+\beta)^2 [2]_q[3]_q} \\
&&\times \bigg[[2]_q[3]_q+2\alpha [2]_q[3]_q +2q[3]_q-2\alpha(1-q)[3]_q -(1-q^2)[2]_q \bigg] x + \frac{q[n+l]_q[n+l-1]_q}{([n+1]_q+\beta)^2[2]_q[3]_q} \\
&&\times \bigg[[2]_q[3]_q-2(1-q)[3]_q+(1-q)^2[2]_q \bigg]x^2 \\
&=& \frac{1}{[3]_q([n+1]_q+\beta)^2} +\frac{2\alpha}{[2]_q([n+1]_q+\beta)^2}+\frac{\alpha^2}{([n+1]_q+\beta)^2}+ \frac{[n+l]_q}{([n+1]_q+\beta)^2 [2]_q[3]_q} \\
&&\times q\bigg((1+q)(1+2q)+(1+q+q^2)(4\alpha+2)\bigg)x + \frac{q[n+l]_q[n+l-1]_q}{([n+1]_q+\beta)^2[2]_q[3]_q} \\
&&\times \bigg((1+q+q^2)(3q-1)+(1+q^2-2q)(1+q)\bigg)x^2 \\
&=&\frac{1}{[3]_q([n+1]_q+\beta)^2} +\frac{2\alpha}{[2]_q([n+1]_q+\beta)^2}+\frac{\alpha^2}{([n+1]_q+\beta)^2} \\
&&+\frac{q[n+l]_q \big((3+4\alpha)+(5+4\alpha)q+4(1+\alpha)q^2\big)}{([n+1]_q+\beta)^2[2]_q[3]_q}x
+ \frac{q^2[n+l]_q[n+l-1]_q(1+q+4q^2)}{([n+1]_q+\beta)^2[2]_q[3]_q}x^2.
\end{eqnarray*}
Hence the lemma.
\newline

Remark 2.1. From the Lemma 2.1, we have
\begin{enumerate}
\item [(i)]$L_{n,l}^{\alpha,\beta}((t-x);q;x)= \bigg(\frac{2q[n+l]_q}{[2]_q ([n+1]_q+\beta}-1\bigg)x+ \frac{1}{[2]_q([n+1]_q+\beta)}+\frac{\alpha}{[n+l]_q}$,
\item [(ii)]$L_{n,l}^{\alpha,\beta}((t-x)^2;q;x)=\frac{\alpha^2}{([n+1]_q+\beta)^2}+\frac{2\alpha}{([n+1]_q+\beta)^2[2]_q} +\frac{1}{([n+1]_q+\beta)^2[3]_q} \\
 ~~~~~~~~~~~~~~~~~~~~~~~~~~+\bigg(\frac{q[n+l]_q \big((3+4\alpha)+(5+4\alpha)q+4(1+\alpha)q^2\big)}{([n+1]_q+\beta)^2[2]_q[3]_q}
 -\frac{2}{([n+1]_q+\beta)^2[2]_q}-\frac{2\alpha}{[n+l]_q}\bigg)x \\
 ~~~~~~~~~~~~~~~~~~~~~~~~~~+\bigg(\frac{q^2[n+l]_q[n+l-1]_q(1_q+4q^2)}{[n+1]_q+\beta)^2[2]_q[3]_q}-\frac{4q[n+l]_q}{[n+1]_q+\beta)[2]_q}+1\bigg)x^2$.
\end{enumerate}

\begin{lemma}
  For any $f\in C(I)$, we have $\Vert L_{n,l}^{\alpha,\beta}(f;q;.)\Vert_{C[0,1]}\leq \Vert f\Vert_{C[0,1+l]}$.
\end{lemma}

\section{Direct Theorems}
\hspace{8mm}In this section, we prove some direct theorems for the operators $L_{n,l}^{\alpha,\beta}(f;q;x)$.
\begin{theorem}
Let $f\in C(I)$ and $0<q_n<1$. Then the sequence of the operators $L_{n,l}^{\alpha,\beta}(f;q_n;.)$ converges uniformly to $f$ on the compact
interval $[0,1]$ if and only if $\lim\limits_{n\rightarrow \infty}q_n =1$.
\end{theorem}
\parindent=0mm\textbf{Proof}.
(Forward) Suppose that $\lim\limits_{n\rightarrow \infty}q_n =1$. Then we shall show that ${L_{n,l}^{\alpha,\beta}(f;q_n;.)}$ converges to $f$
uniformly on $[0,1]$. Note that for $0<q_n<1$ and $q_n\rightarrow\infty$ for $n\rightarrow\infty$, we get $[n+1]_{q_n}\rightarrow\infty$
as $n\rightarrow\infty$. Now it is easily seen that $\frac{[n+l]_{q_n}}{[n+1]_{q_n}+\beta}=1+q_n^n \frac{([l]_{q_n}-1)}{[n+1]_{q_n}}
-\frac{\beta}{[n+1]_{q_n}+\beta}$. So when $n\rightarrow \infty$, $\frac{[n+l]_{q_n}}{[n+1]_{q_n}+\beta} \rightarrow 1$ and
$\frac{[n+l]_{q_n}}{([n+1]_{q_n}+\beta)^2} \rightarrow 0$. Using this and the Lemma 2.1, we find that $L_{n,l}^{\alpha,\beta}(1;q_n;x)\rightarrow 1$,
$L_{n,l}^{\alpha,\beta}(t;q_n;x)\rightarrow x$ and $L_{n,l}^{\alpha,\beta}(t^2;q_n;x)\rightarrow x^2$ uniformly on the compact set $[0, 1]$
as ${n\rightarrow \infty}$. Therefore, the Korovkin's theorem proves that the sequence $L_{n,l}^{\alpha,\beta}(f;q_n;.)$ converges uniformly to $f$ on $[0, 1]$.
\newline
We shall prove the converse by the method of contradiction. Suppose that the sequence $(q_n)$ does not converge to $1$. Then there must exist a subsequence
$(q_{n_i})$ of the sequence $(q_n)$ such that $q_{n_i}\in (0, 1)$, $q_{n_i} \rightarrow \delta \in [0,1)$ as $i \rightarrow \infty$.
 Then $\frac{1}{[n_i+p]_{q_{n_i}}} = \frac{1-q_{n_i}}{1-(q_{n_i})^{n_i+p}}\rightarrow 1-\delta$ as $i\rightarrow \infty$ because $(q_{n_i})^{n_i} \rightarrow 0$
 as $i \rightarrow \infty$. Now if we choose $n= n_i, q = q_{n_i}$ in $L_{n,l}^{\alpha,\beta}(t;q;x)$ form the Lemma 2.1, then we get
 $L_{n,l}^{\alpha,\beta}(t;q;x) = \frac{2\delta}{(1+\delta)(1+\beta(1-\delta))}x +\frac{1-\delta}{1+\delta} \frac{1}{(1+\beta(1-\delta))}+ \alpha (1-\delta)$,
 which is different from $x$ when $i \rightarrow \infty$, which contradicts our supposition. Hence $\lim\limits_{n\rightarrow \infty}q_n =1$. Hence the theorem is completely proved.
 \newline

Now we define the following:\\
 Let $f\in C(I)$, $\delta>0$ and $W^2 = \left\{h: h', h'' \in C(I)\right\}$, then the Peetre's K-functional is defined by
\begin{equation*}
K_2(f,\delta)= \inf\limits_{h\in W^2}\{\Vert f-h\Vert+\delta \Vert h''\Vert\},
\end{equation*}
By DeVore and Lorentz theorem [?]
 there exists a constant $C>0$ such that
\begin{equation}
K_2(f,\delta)\leq C \omega_2(f, \sqrt{\delta})
\end{equation}
where $\omega_2 (f, \sqrt{\delta})$, the second order modulus of continuity of $f\in C(I)$, is defined as
\begin{equation*}
\omega_2(f, \sqrt{\delta})=\sup\limits_{0<p<\delta^\frac{1}{2}}\sup\limits_{x\in I}|f(x+2p)-2f(x+p)+f(x)|.
\end{equation*}
Also by $\omega (f, \delta)$, we denote the first order modulus of continuity of $f\in C(I)$ defined as
\begin{equation*}
\omega (f, \delta) = \sup\limits_{0<p<\delta}\sup\limits_{x\in I}|f(x+p)-f(x)|
\end{equation*}
Next we prove the following theorem.
\begin{theorem}\label{theorem3.2}
Let $L_{n,l}^{\alpha, \beta} (f; q; x)$ be the sequence of positive linear operators defined by (2.1) and $f\in C(I)$. Let $(q_n)$ be the sequence
with $0<q_n<1$ and $q_n\rightarrow 1$ as $n\rightarrow \infty$. Then there exists a constant $c>0$ indepenent of $n$ and $x$ such that
\begin{equation}
\bigl{|}L_{n,l}^{\alpha, \beta} (f; q; x)-f(x)\bigl{|}\leq C \omega_2\left(f,\sqrt {\phi_{n,l}^{\alpha, \beta}(q_n; x)}\right)
+\omega \left(f, \frac{\alpha}{[n+l]_{q_n}}+\frac{1}{[2]_{q_n}([n+1]_{q_n}+\beta)}+\frac{2q_n[n+l]_{q_n}}{[2]_{q_n}([n+1]_{q_n}+\beta)}x-x\right)
\end{equation}
where $\phi_{n,l}^{\alpha, \beta}(q_n; x) = L_{n,l}^{\alpha, \beta} ((t-x)^2; q_n; x)
+\left(\frac{\alpha}{[n+l]_{q_n}}+\frac{1}{[2]_{q_n}([n+1]_{q_n}+\beta)}+\frac{2q_n[n+l]_{q_n}}{[2]_{q_n}([n+1]_{q_n}+\beta)}x-x\right)^2$ and $x\in [0, 1]$.
\end{theorem}
\parindent=0mm\textbf{Proof}. Let us define the following operators
\begin{equation}
\bar{L}_{n,l}^{\alpha, \beta} (f; q_n; x) = L_{n,l}^{\alpha, \beta} (f; q_n; x)+ f(x)- f\left(\frac{\alpha}{[n+l]_{q_n}}+\frac{1}{[2]_{q_n}([n+1]_{q_n}+\beta)}
+\frac{2q_n[n+l]_{q_n}}{[2]_{q_n}([n+1]_{q_n}+\beta)}x\right)
\end{equation}
In the light of the Lemma 2.1, it is easily seen that $\bar{L}_{n,l}^{\alpha, \beta} (1; q_n; x) =1$ and $\bar{L}_{n,l}^{\alpha, \beta} (t; q_n; x) =x$.
Now from the Taylor's formula, for $g\in W^2$, we can write
\begin{equation*}
g(t)-g(x) = (t-x)g'(x)+ ~\int_{x}^{t}(t-u)g''(u)du
\end{equation*}
Applying the operators $\bar{L}_{n,l}^{\alpha, \beta}$ to both sides of the above equation  we get
\begin{eqnarray*}
\bar{L}_{n,l}^{\alpha, \beta}(g; q_n; x)- g(x)& =& g'(x)\bar{L}_{n,l}^{\alpha, \beta}((t-x); q_n; x)
+ \bar{L}_{n,l}^{\alpha, \beta} \left( ~\int_{x}^{t}(t-u)g''(u)du\right) \\
&=&\bar{L}_{n,l}^{\alpha, \beta} \left( ~\int_{x}^{t}(t-u)g''(u)du, q_n; x\right)\\
&=&L_{n,l}^{\alpha, \beta} \left( ~\int_{x}^{t}(t-u)g''(u)du, q_n; x\right)\\
&&- ~\int_{x}^{\frac{\alpha}{[n+l]_{q_n}}+\frac{1}{[2]_{q_n}([n+1]_{q_n}+\beta)}
+\frac{2q_n[n+l]_{q_n}}{[2]_{q_n}([n+1]_{q_n}+\beta)}x}
\bigg(\frac{\alpha}{[n+l]_{q_n}}+\frac{1}{[2]_{q_n}([n+1]_{q_n}+\beta)}\\
&&+\frac{2q_n[n+l]_{q_n}}{[2]_{q_n}([n+1]_{q_n}+\beta)}x -u\bigg) g''(u)du
\end{eqnarray*}
Therefore, we will have
\begin{eqnarray*}
\vert \bar{L}_{n,l}^{\alpha, \beta}(g; q_n; x)- g(x)\vert &\leq& L_{n,l}^{\alpha, \beta}((t-x)^2; q_n; x) \Vert g''\Vert_{C[0, 1+l]}
+\bigg(\frac{\alpha}{[n+l]_{q_n}}+\frac{1}{[2]_{q_n}([n+1]_{q_n}+\beta)}\\
&&+\frac{2q_n[n+l]_{q_n}}{[2]_{q_n}([n+1]_{q_n}+\beta)}x-x\bigg)^2\Vert g''\Vert_{C[0, 1+l]}\\
&=& \phi_{n,l}^{\alpha, \beta}(q_n; x)\Vert g''\Vert_{C[0, 1+l]}
\end{eqnarray*}
In view of (3.2), we obtain
\begin{eqnarray*}
 \vert L_{n,l}^{\alpha, \beta}(f; q_n; x)-f(x)\vert &\leq& \vert \bar{L}_{n,l}^{\alpha, \beta}(f-g; q_n; x)- g(x)\vert +
 \vert \bar{L}_{n,l}^{\alpha, \beta}(g; q_n; x)- g(x)\vert \\
&+&\bigg\vert f\left(\frac{\alpha}{[n+l]_{q_n}}+\frac{1}{[2]_{q_n}([n+1]_{q_n}+\beta)}
+\frac{2q_n[n+l]_{q_n}}{[2]_{q_n}([n+1]_{q_n}+\beta)}x\right)-f(x)\bigg \vert
\end{eqnarray*}
Now we have
\begin{equation*}
\Vert \bar{L}_{n,l}^{\alpha, \beta}(f; q_n; x)\Vert \leq 3 \Vert f\Vert_{C[0, 1+l]},
\end{equation*}
 by the Lemma 2.2, so we have \\
$\vert L_{n,l}^{\alpha, \beta}(f; q_n; x)-f(x)\vert \leq 4\Vert f-g\Vert_{C[0, 1+l]}
+\phi_{n,l}^{\alpha, \beta}(q_n; x)\Vert g''\Vert_{C[0, 1+l]} + \omega \bigg( f, \bigg\vert \frac{1}{([n+1]_{q_n}+\beta)[2]_{q_n}}\\
~~~~~~~~~~~~~~~~~~~~~~~~~~~~~~~~~~~+\frac{[n+l]_{q_n}}{([n+1]_{q_n}+\beta)} \frac{2q_n}{[2]_{q_n}}x -x\bigg\vert \bigg)$.\\
On taking the infimum of the right hand side running over all $g\in W^2$ and using the definition of the Peetre's functional, we get
\begin{equation*}
\vert L_{n,l}^{\alpha, \beta}(f; q_n; x)-f(x)\vert \leq 4 K_2(f,\phi_{n,l}^{\alpha, \beta}(q_n; x))
+ \omega \bigg(f,\frac{1}{([n+1]_{q_n}+\beta)[2]_{q_n}}
+\frac{[n+l]_{q_n}}{([n+1]_{q_n}+\beta)} \frac{2q_n}{[2]_{q_n}}x -x\bigg\vert \bigg)
\end{equation*}
Now in view of (3.1), we obtain
\begin{equation*}
\vert L_{n,l}^{\alpha, \beta}(f; q_n; x)-f(x)\vert \leq C\omega_2\bigg(f,\sqrt{\phi_{n,l}^{\alpha, \beta}(q_n; x)}\bigg)
+ \omega \bigg(f,\bigg\vert\frac{1}{([n+1]_{q_n}+\beta)[2]_{q_n}}
+\frac{[n+l]_{q_n}}{([n+1]_{q_n}+\beta)} \frac{2q_n}{[2]_{q_n}}x -x\bigg\vert \bigg),
\end{equation*}
and this completes the proof of the theorem.\\

Now we shall obtain an estimate of the rate of convergence for the operators defined in (2.1) using the Lipschitz-type maximal function defined as follows [?]\\
For $x\in[0,1]$ and $\xi \in (0,1]$, the Lipschitz-type maximal function is defined as
\begin{equation}
\tilde{\omega}_\xi(f, x) = \sup\limits_{t\neq x, t\in[0, 1+l]} \frac{\vert f(t)-f(x)\vert}{\vert t-x\vert^{\xi}}
\end{equation}
Now we prove the following theorem
\begin{theorem}\label{theorem3.3}
Let $f\in C(I), 0<\xi\leq1$ and $q_n\in (0,1)$ such that $q_n \rightarrow 1$ as $n\rightarrow \infty$. Then for every $x\in [0,1]$, we have
\begin{equation*}
 \vert L_{n,l}^{\alpha, \beta}(f; q_n; x)-f(x)\vert \leq \tilde{\omega}_\xi(f, x)(\gamma_{n,l}(q_n; x))^\frac{\xi}{2},
\end{equation*}
where $\gamma_{n,l}(q_n; x) = L_{n,l}^{\alpha, \beta}((t-x)^2; q_n; x)$
\end{theorem}
\parindent=0mm\textbf{Proof}. In the light of the Lemma (2.1), we have
\begin{equation}
\vert L_{n,l}^{\alpha, \beta}(f; q_n; x)-f(x)\vert \leq L_{n,l}^{\alpha, \beta}(\vert f(t)-f(x)\vert; q_n; x) \leq \tilde{\omega}_\xi(f, x)
L_{n,l}^{\alpha, \beta}(\vert t-x\vert^{\xi}; q_n; x)
\end{equation}
and in view of (3.4), we have
\begin{equation}
\vert f(t)-f(x)\vert  \leq \tilde{\omega}_\xi(f, x)\vert t-x \vert^\xi.
\end{equation}
When we use the H\"{o}lder's inequality with $\tilde{p}= \frac{2}{\xi}$ and $\tilde{q}= \frac{2}{2-\xi}$, we obtain
\begin{equation*}
\vert L_{n,l}^{\alpha, \beta}(f; q_n; x)-f(x)\vert \leq \tilde{\omega}_\xi(f, x)
L_{n,l}^{\alpha, \beta}(\vert t-x\vert^{2}; q_n; x)^{\frac{\xi}{2}} = \tilde{\omega}_\xi(f, x)(\gamma_{n,l}(q_n; x))^\frac{\xi}{2},
\end{equation*}
and hence the theorem.\\

To prove the next theorem we consider the following Lipschitz-type space of functions [?]:
\begin{equation*}
\tilde{Lip}_M(s)= \left\{f\in C(I): \vert f(t)-f(x)\vert \leq M\frac{\vert t-x\vert^s}{(t+x)^{\frac{s}{2}}}\right\},
\end{equation*}
where $M$ is a positive constant and $0<s\leq1$. \\
Now we have the following theorem.
\begin{theorem}\label{theorem3.4}
Let $f\in \tilde{Lip}_M(s), s\in (0, 1]$ and $q_n\in (0,1)$ such that $q_n \rightarrow 1$ as $n\rightarrow \infty$. Then for each $x\in(0, 1]$, we have
\begin{equation*}
\vert L_{n,l}^{\alpha, \beta}(f; q_n; x)-f(x)\vert \leq M \left(\frac{\gamma_{n,l}(q_n; x)}{x}\right)^{\frac{s}{2}}
\end{equation*}
where $\gamma_{n,l}(q_n; x) =  L_{n,l}^{\alpha, \beta}((t-x)^2; q_n; x)$.
\end{theorem}
\parindent=0mm\textbf{Proof}.
Firstly we will prove the result for $s=1$. In fact we have, for $f\in \tilde{Lip}_M(1)$,
\begin{eqnarray*}
\vert L_{n,l}^{\alpha, \beta}(f; q_n; x)-f(x)\vert &\leq&([n+1]_{q_n}+\beta)\sum\limits_{k=0}^{n+l}b_{n,l}^k(q_n;x)q_n^{-k}
\int_{\frac{[k]_{q_n}+\alpha}{[n+1]_{q_n}+\beta}}^{\frac{[k+1]_{q_n}+\alpha}{[n+1]_{q_n}+\beta}}\vert f(t)-f(x)\vert d_{q_n}^Rt \\
&\leq& M ([n+1]_{q_n}+\beta)\sum\limits_{k=0}^{n+l}b_{n,l}^k(q_n;x)q_n^{-k}
\int_{\frac{[k]_{q_n}+\alpha}{[n+1]_{q_n}+\beta}}^{\frac{[k+1]_{q_n}+\alpha}{[n+1]_{q_n}+\beta}}\frac{\vert t-x\vert}{\sqrt{t+x}} d_{q_n}^Rt
\end{eqnarray*}
Applying the Cauchy-Schwarz inequality and the fact $\frac{1}{\sqrt{t+x}}\leq \frac{1}{\sqrt {x}}$, we get \\
\begin{eqnarray*}
\vert L_{n,l}^{\alpha, \beta}(f; q_n; x)-f(x)\vert &\leq& \frac{M}{\sqrt{x}}([n+1]_{q_n}+\beta)\sum\limits_{k=0}^{n+l}b_{n,l}^k(q_n;x)q_n^{-k}
\int_{\frac{[k]_{q_n}+\alpha}{[n+1]_{q_n}+\beta}}^{\frac{[k+1]_{q_n}+\alpha}{[n+1]_{q_n}+\beta}}\vert t-x\vert d_{q_n}^Rt \\
&=& \frac{M}{\sqrt{x}} L_{n,l}^{\alpha, \beta}(\vert t-x \vert; q_n; x) \\
&\leq& M \left(\frac{\gamma_{n,l}(q_n; x)}{x}\right)^{\frac{1}{2}}.
\end{eqnarray*}
Thus the result is established for $s=1$. Next we prove the result for $0<s<1$. On using the H\"{o}lder's inequality twice for $\tilde{p}=\frac{1}{s}$
and $\tilde{q} = \frac{1}{1-s}$, we get
\begin{eqnarray*}
\vert L_{n,l}^{\alpha, \beta}(f; q_n; x)-f(x)\vert &\leq& ([n+1]_{q_n}+\beta)\sum\limits_{k=0}^{n+l}b_{n,l}^k(q_n;x)q_n^{-k}
\int_{\frac{[k]_{q_n}+\alpha}{[n+1]_{q_n}+\beta}}^{\frac{[k+1]_{q_n}+\alpha}{[n+1]_{q_n}+\beta}}\vert f(t)-f(x)\vert d_{q_n}^Rt \\
&\leq & \Bigg\{\sum\limits_{k=0}^{n+l}b_{n,l}^k(q_n;x)\left(([n+1]_{q_n}+\beta)q_n^{-k}
\int_{\frac{[k]_{q_n}+\alpha}{[n+1]_{q_n}+\beta}}^{\frac{[k+1]_{q_n}+\alpha}{[n+1]_{q_n}+\beta}}\vert f(t)-f(x)\vert d_{q_n}^Rt \right)^{\frac{1}{s}}\Bigg\}^s \\
&\leq & \Bigg\{([n+1]_{q_n}+\beta)\sum\limits_{k=0}^{n+l}b_{n,l}^k(q_n;x)q_n^{-k}
\int_{\frac{[k]_{q_n}+\alpha}{[n+1]_{q_n}+\beta}}^{\frac{[k+1]_{q_n}+\alpha}{[n+1]_{q_n}+\beta}}\vert f(t)-f(x)\vert^{\frac{1}{s}} d_{q_n}^Rt\Bigg\}^s.
\end{eqnarray*}
Now as $f\in \tilde{Lip}_M(s)$, we obtain
\begin{eqnarray*}
\vert L_{n,l}^{\alpha, \beta}(f; q_n; x)-f(x)\vert &\leq& M \Bigg\{([n+1]_{q_n}+\beta)\sum\limits_{k=0}^{n+l}b_{n,l}^k(q_n;x)q_n^{-k}
\int_{\frac{[k]_{q_n}+\alpha}{[n+1]_{q_n}+\beta}}^{\frac{[k+1]_{q_n}+\alpha}{[n+1]_{q_n}+\beta}}\frac{\vert t-x\vert}{\sqrt{t+x}} d_{q_n}^Rt\Bigg\}^s \\
&\leq& \frac{M}{x^{\frac{s}{2}}} \Bigg\{([n+1]_{q_n}+\beta)\sum\limits_{k=0}^{n+l}b_{n,l}^k(q_n;x)q_n^{-k}
\int_{\frac{[k]_{q_n}+\alpha}{[n+1]_{q_n}+\beta}}^{\frac{[k+1]_{q_n}+\alpha}{[n+1]_{q_n}+\beta}}\vert t-x\vert d_{q_n}^Rt\Bigg\}^s \\
&=& \frac{M}{x^{\frac{s}{2}}}(L_{n,l}^{\alpha, \beta}(\vert t-x \vert; q_n; x))^s \\
&\leq& M \left(\frac{\gamma_{n,l}(q_n; x)}{x}\right)^{\frac{s}{2}}.
\end{eqnarray*}
This completes the proof of the theorem.

\section{Statistical Convergence}
In this section we study the $A$-statistical convergence of the operators defined in (2.1) through the Korovkin-type statistical approximation theorem.\\
Let $A=(a_{nk})$ be a non-negative infinite summability matrix. Then for a sequence $x=(x_k)$, we define the $A$-transform of $x$ as
$(Ax)_n = \sum_{k=1}^{\infty}$ whenever the series converges for each $n$. We denote it by $Ax = \left((Ax)_n\right)$. $A$ is said to be regular if $\lim_n (Ax)_n
=L$, whenever $\lim_n (x)_n = L$. The sequence $x=(x)_n$ is said to be $A$-statistically convergent to the limit $L$, denoted by $st_A-\lim_n x_n = L$, if for
each $\epsilon >0,  \lim_n \sum_{k:\vert x_k-L\vert\geqslant\epsilon} a_{nk} = 0$. If we take the matrix $A$ to be the $Ces\grave{a}ro$ matrix $C$ of order one,
then the notion of the $A$-statistically convergence is same as the statistical convergence. \\
Now to prove a theorem, we take a sequence $(q)_n$ such that $q_n \in (0,1)$ satisfying the following:
$st_A-\lim_n x_n = 1$, $st_A-\lim_n (q_n)^n = a \in (0,1)$, $st_A-\lim_n \frac{1}{[n]_{q_n}} = 0$.
\begin{theorem} \label{4.1}
Let $A=(a_{nk})$ be a non-negative regular summability matrix and $(q_n)$ be a sequence satisfying the above conditions. Then for any $f\in C(I)$, we have
\begin{equation*}
st_A-\lim_n \Vert L_{n,l}^{\alpha,\beta}(f;q;.)-f\Vert_{C[0,1]} = 0.
\end{equation*}
\end{theorem}
\parindent=0mm\textbf{Proof}. Let $e_i(x) = x^i$, where $x\in[0,1]$, $i =0, 1, 2$. Then from the Lemma 2.1, we have
\begin{equation}
st_A-\lim_n \Vert L_{n,l}^{\alpha,\beta}(e_0;q;.)-e_0\Vert_{C[0,1]} = 0.
\end{equation}
Next, again from the Lemma 2.1, we have
\begin{equation*}
\lim_n \Vert L_{n,l}^{\alpha,\beta}(e_1;q_n;.)-e_1\Vert_{C[0,1]}\leq \left\vert\frac{\alpha}{[n+l]_q}+\frac{1}{([n+1]_q+\beta)[2]_q}\right\vert
+\left\vert\frac{2q[n+l]_q}{([n+1]_q+\beta)[2]_q}-1\right\vert.
\end{equation*}
Now since $st_A-\lim_n q_n = 1$, $st_A-\lim_n (q_n)^n = a \in (0,1)$ and $st_A-\lim_n \frac{1}{[n]_{q_n}} = 0$, we have
\begin{equation*}
st_A-\lim_n \left(\frac{\alpha}{[n+l]_q}+\frac{1}{([n+1]_q+\beta)[2]_q}\right)=0
\end{equation*}
and
\begin{equation*}
st_A-\lim_n \left(\frac{2q_n}{[2]_{q_n}}\frac{1-q_n^{n+l}}{1-q_n^{n+1}+\beta (1-q_n)}-1\right)=0
\end{equation*}
Now, for a given $\epsilon>0$, let us define the following sets:
\begin{equation*}
U= \left\{n\in N:\Vert L_{n,l}^{\alpha,\beta}(e_1;q_n;.)-e_1\Vert_{C[0,1]}\geq \epsilon \right\},
\end{equation*}
\begin{equation*}
U_1= \left\{n\in N:\frac{\alpha}{[n+l]_{q_n}}+\frac{1}{([n+1]_{q_n}+\beta)[2]_{q_n}}\geq \frac{\epsilon}{2} \right\},
\end{equation*}
and
\begin{equation*}
U_2= \left\{n\in N:\frac{2q_n}{[2]_{q_n}}\frac{1-q_n^{n+l}}{1-q_n^{n+1}+\beta (1-q_n)}-1 \geq \frac{\epsilon}{2}\right\}.
\end{equation*}
The containment $U\subseteq U_1\cup U_2$ is obvious which in turn implies that $\sum_{n\in U}a_{nk}\leq\sum_{n\in U_1}a_{nk}+\sum_{n\in U_2}a_{nk}$, and
hence we have
\begin{equation}
st_A-\lim_n \Vert L_{n,l}^{\alpha,\beta}(e_1;q_n;.)-e_1\Vert_{C[0,1]}.
\end{equation}
Further, using the Lemma 2.1, we have
\begin{eqnarray*}
\Vert L_{n,l}^{\alpha,\beta}(e_2;q_n;.)-e_2\Vert_{C[0,1]}\leq\left\vert\frac{1}{([n+1]_{q_n}+\beta)^2[3]_{q_n}}+\frac{2\alpha}{([n+1]_q+\beta)^2[2]_{q_n}}
+\frac{\alpha^2}{([n+1]_{q_n}+\beta)^2}\right\vert \\
+\left\vert\frac{q_n}{[2]_{q_n}[3]_{q_n}}\frac{((3+4\alpha)+(5+4\alpha)q_n+4(1+\alpha)q_n^2)(1-q_n^{n+l})(1-q_n))}{(1-q_n^{n+1}+\beta(1-q_n))^2}\right\vert \\
+ \left\vert \frac{q_n^2(1+q_n+4q_n^2)}{[2]_{q_n}[3]_{q_n}} \frac{1-q_n^{n+l}}{(1-q_n^{n+1}+\beta (1-q_n))}\frac{1-q_n^{n+l-1}}{(1-q_n^{n+1}
+\beta (1-q_n))}-1\right\vert.
\end{eqnarray*}
Keeping in view $st_A-\lim_n q_n = 1$, $st_A-\lim_n (q_n)^n = a \in (0,1)$ and $st_A-\lim_n \frac{1}{[n]_{q_n}} = 0$, we obtain
\begin{equation*}
st_A-\lim_n \left(\frac{1}{([n+1]_{q_n}+\beta)^2[3]_{q_n}}+\frac{2\alpha}{([n+1]_q+\beta)^2[2]_{q_n}}
+\frac{\alpha^2}{([n+1]_{q_n}+\beta)^2}\right)=0,
\end{equation*}
\begin{equation*}
st_A-\lim_n \left(\frac{q_n}{[2]_{q_n}[3]_{q_n}}\frac{((3+4\alpha)+(5+4\alpha)q_n+4(1+\alpha)q_n^2)(1-q_n^{n+l})(1-q_n))}{(1-q_n^{n+1}
+\beta(1-q_n))^2}\right)=0,
\end{equation*}
and
\begin{equation*}
st_A-\lim_n \left(\frac{q_n^2(1+q_n+4q_n^2)}{[2]_{q_n}[3]_{q_n}} \frac{1-q_n^{n+l}}{(1-q_n^{n+1}+\beta (1-q_n))}\frac{1-q_n^{n+l-1}}{(1-q_n^{n+1}
+\beta (1-q_n))}-1\right)=0.
\end{equation*}
Now for each $\epsilon >0$, we define the following sets:
\begin{equation*}
V= \left\{n\in N:\Vert L_{n,l}^{\alpha,\beta}(e_2;q_n;.)-e_2\Vert_{C[0,1]}\geq \epsilon \right\},
\end{equation*}
\begin{equation*}
V_1= \left\{n\in N:\frac{1}{([n+1]_{q_n}+\beta)^2[3]_{q_n}}+\frac{2\alpha}{([n+1]_{q_n}+\beta)^2[2]_{q_n}}
+\frac{\alpha^2}{([n+1]_{q_n}+\beta)^2}\geq\frac{\epsilon}{3}\right\},
\end{equation*}
\begin{equation*}
V_2= \left\{n\in N:\frac{q_n}{[2]_{q_n}[3]_{q_n}}\frac{((3+4\alpha)+(5+4\alpha)q_n+4(1+\alpha)q_n^2)(1-q_n^{n+l})(1-q_n))}{(1-q_n^{n+1}
+\beta(1-q_n))^2}\geq \frac{\epsilon}{3}\right\},
\end{equation*}
and
\begin{equation*}
V_3= \left\{n\in N:\frac{q_n^2(1+q_n+4q_n^2)}{[2]_{q_n}[3]_{q_n}} \frac{1-q_n^{n+l}}{(1-q_n^{n+1}+\beta (1-q_n))}\frac{1-q_n^{n+l-1}}{(1-q_n^{n+1}
+\beta (1-q_n))}-1\geq\frac{\epsilon}{3}\right\},
\end{equation*}
It is obvious that $V\subseteq V_1\cup V_2\cup V_3$, which in turn implies that $\sum_{n\in V}a_{nk}\leq\sum_{n\in V_1}a_{nk}+\sum_{n\in V_2}a_{nk}
+\sum_{n\in V_3}a_{nk}$. Therefore, we get
\begin{equation}
st_A-\lim_n \Vert L_{n,l}^{\alpha,\beta}(e_2;q_n;.)-e_2\Vert_{C[0,1]}=0.
\end{equation}
Now on combining (4.1)-(4.3), the theorem follows from the Korovkin-type statistical approximation theorem as proved in [?]. Hence the proof is complete.

\section{Construction of the Bivariate Operators}
In what follows we construct the bivariate extension of the operators defined by (2.1). \\
Let $I_1 = [0,1+l_1]$ and $I_2 = [0, 1+l_2]$. We consider $C(I_1\times I_2)$, the space of all real valued continuous functions defined on $I_1 \times I_2$ equipped
with the following norm
\begin{equation*}
\Vert f\Vert_{C(I_1 \times I_2)} = \sup\limits_{(x, y)\in I_1\times I_2} \vert f(x, y)\vert.
\end{equation*}
We define the bivariate generaliation of the operators in (2.1) as follows
\begin{eqnarray*}
L_{n_1,n_2; l_1, l_2}^{\alpha_1,\alpha_2; \beta_1, \beta_2}(f(t, s);q_1, q_2;x, y)&=&([n_1+1]_{q_1}+\beta_1)([n_2+1]_{q_2}+\beta_2)
\sum\limits_{k_1=0}^{n_1+l_1} \sum\limits_{k_2=0}^{n_2+l_2}q_1^{-k_1} q_2^{-k_2} \\
&&\times b_{n_1,n_2; l_1, l_2}^{k_1,k_2}(q_1, q_2;x, y)\int_{\frac{[k_1]_{q_1}+\alpha_1}
{[n_1+1]_{q_1}+\beta_1}}^{\frac{[k_1+1]_{q_1}+\alpha_1}{[n_1+1]_{q_1}+\beta}}
\int_{\frac{[k_2]_{q_2}+\alpha_2}{[n+1]_{q_2}+\beta_2}}^{\frac{[k+1]_{q_2}+\alpha_2}{[n+1]_{q_2}+\beta_2}}
f(t, s)d_{q_1}^Rt d_{q_2}^Rt,
\end{eqnarray*}
where
\begin{eqnarray*}
b_{n_1,n_2; l_1, l_2}^{k_1,k_2}(q_1, q_2;x, y) = \binom{n_1+l_1}{k_1}_{q_1} \binom{n_2+l_2}{k_2}_{q_2} x^{k_1} y^{k_2}(1-x)_{q_1}^{n_1+l_1-k_1}
(1-x)_{q_2}^{n_2+l_2-k_2},
\end{eqnarray*}
$f\in C(I_1 \times I_2)$, $0< q_1, q_2<1$, $(x, y) \in [0,1] \times [0,1] = J^2$ and $\alpha_1, \alpha_2, \beta_1, \beta_2$ are such
that ~$0<\beta_1\leq\alpha_1; 0<\beta_2\leq\alpha_2$. \\
Now we prove a lemma concerning the above operators.
\begin{lemma}
Let $(t, s) \in (I_1\times I_2)$, $(i, j) \in N^0\times N^0$ with $i+j \leq2$, and $t^js^j$ by $e_{ij}(t,s)$ be the two dimensional test functions.
Then the following equalities hold for the bivariate operators of (?):
\end{lemma}
\begin{enumerate}
\item [(i)] $L_{n_1,n_2; l_1, l_2}^{\alpha_1,\alpha_2; \beta_1, \beta_2}(e_{00};q_1, q_2;x, y)=1$,
\item [(ii)]$L_{n_1,n_2; l_1, l_2}^{\alpha_1,\alpha_2; \beta_1, \beta_2}(e_{10};q_1, q_2;x, y)=\frac{\alpha_1}{[n_1+l_1]_{q_1}}+\frac{1}{([n_1+1]_{q_1}
+\beta_1)[2]_{q_1}}+\frac{2q_1[n_1+l_1]_{q_1}}{([n_1+1]_{q_1}+\beta_1)[2]_{q_1}}x$,
\item [(iii)]$L_{n_1,n_2; l_1, l_2}^{\alpha_1,\alpha_2; \beta_1, \beta_2}(e_{01};q_1, q_2;x, y)=\frac{\alpha_2}{[n_2+l_2]_{q_2}}+\frac{1}{([n_2+1]_{q_2}
+\beta_2)[2]_{q_2}}+\frac{2q_2[n_2+l_2]_{q_2}}{([n_2+1]_{q_2}+\beta_2)[2]_{q_2}}y$,
\item [(iv)]$L_{n_1,n_2; l_1, l_2}^{\alpha_1,\alpha_2; \beta_1, \beta_2}(e_{20};q_1, q_2;x, y)=\frac{1}{([n_1+1]_{q_1}+\beta_1)^2[3]_{q_1}}
+\frac{2\alpha_1}{([n_1+1]_{q_1}+\beta_1)^2[2]_{q_1}}+\frac{\alpha_1^2}{([n_1+1]_{q_1}+\beta_1)^2} \\
    ~~~~~~~~~~~~~~~~~~~~~~~~~~~~~~~~~~~~~~~+\frac{{q_1}[n_1+l_1]_{q_1}((3+4\alpha_1)+(5+4\alpha_1){q_1}+4(1+\alpha_1){q_1}^2)}{([n_1+1]_{q_1}
    +\beta_1)^2[2]_{q_1}[3]_{q_1}}x
    +\frac{{q_1}^2[n_1+l_1]_{q_1}[n_1+l_1-1]_{q_1}(1+{q_1}+4{q_1}^2)}{([n_1+1]_{q_1}+\beta_1)^2[2]_{q_1}[3]_{q_1}}x^2$,
\item [(v)]$L_{n_1,n_2; l_1, l_2}^{\alpha_1,\alpha_2; \beta_1, \beta_2}(e_{02};q_1, q_2;x, y)=\frac{1}{([n_2+1]_{q_2}+\beta_2)^2[3]_{q_2}}
+\frac{2\alpha_2}{([n_2+1]_{q_2}+\beta_2)^2[2]_{q_2}}+\frac{\alpha_2^2}{([n_2+1]_{q_2}+\beta_2)^2} \\
    ~~~~~~~~~~~~~~~~~~~~~~~~~~~~~~~~~~~~~~~+\frac{{q_2}[n_2+l_2]_{q_2}((3+4\alpha_2)+(5+4\alpha_2){q_2}+4(1+\alpha_2){q_2}^2)}{([n_2+1]_{q_2}
    +\beta_1)^2[2]_{q_2}[3]_{q_2}}y
    +\frac{{q_2}^2[n_2+l_2]_{q_2}[n_2+l_2-1]_{q_2}(1+{q_2}+4{q_2}^2)}{([n_2+1]_{q_2}+\beta_2)^2[2]_{q_2}[3]_{q_2}}y^2$.
\end{enumerate}
\parindent=0mm\textbf{Proof}. In the light of the Lemma 2.1 and noting that
\begin{equation*}
L_{n_1,n_2; l_1, l_2}^{\alpha_1,\alpha_2; \beta_1, \beta_2}(t^i s^j;q_1, q_2;x, y)
= L_{n_1, l_1}^{\alpha_1, \beta_1}(t^i ;q_1,x)\times L_{n_2, l_2}^{\alpha_2, \beta_2}(s^j ;q_2,y), ~~for~0\leq i, j\leq2.
\end{equation*}
the proof is plain and straightforward. So we omit the details. \\
{} Now To establish the next theorem we first define the following: \\
Let $f\in C(I_1\times I_2)$ and $\delta_1, \delta_2>0$. Then the first order complete modulus of continuity for the bivariate case,
denoted by $\omega (f; \delta_1, \delta_2)$ is defined as follows:
\begin{equation*}
\omega (f, \delta_1, \delta_2) = sup \big\{\vert f(t, s)-f(x, y)\vert:\vert t-x\vert\leq\delta_1, \vert s-y\vert\leq\delta_2\big\}.
\end{equation*}
Two chief properties of $\omega (f, \delta_1, \delta_2)$ are as follows:
\begin{enumerate}
  \item [(i)] $\omega (f, \delta_1, \delta_2) \rightarrow 0 ~as ~\delta_1 \rightarrow 0~ and~ \delta_2 \rightarrow 0$, and
  \item [(ii)] $\vert f(t,s) - f(x, y)\vert \leq \omega(f,\delta_1,\delta_2) \left( 1+\frac{\vert t-x\vert}{\delta_1}\right)
  \left(1+\frac{\vert s-y\vert}{\delta_2}\right)$.
\end{enumerate}
Next we state and prove a theorem regarding the rate of convergence of the bivariate operators. To do it we consider a sequence $(q_{n_i})$
with $q_{n_i}\in (0,1)$ such that $q_{n_i} \rightarrow 1$ and $q_{n_i}^{n_i} \rightarrow a_i, (0\leq a_i<1)$ as $n_i \rightarrow \infty$ for $i = 1, 2$.
Also we denote $L_{n_1, l_1}^{\alpha_1, \beta_1}(t-x)^2 ;q_{n_1},x)$ and $L_{n_2, l_2}^{\alpha_2, \beta_2}(s-y)^2 ;q_{n_2},y)$  by $\delta_{n_1}(x)$ and
$\delta_{n_2}(x)$ respectively. Now we have the follwoing theorem:
\begin{theorem}
Let $f\in C(I_1 \times I_2)$. Then for all $(x,y) \in J^2$, we have
\begin{equation*}
\vert L_{n_1,n_2; l_1, l_2}^{\alpha_1,\alpha_2;\beta_1,\beta_2}(f;q_{n_1}, q_{n_2};x, y)-f(x, y) \vert \leq 4 \omega(f, (\delta_{n_1}(x))^{\frac{1}{2}}
(\delta_{n_2}(y))^{\frac{1}{2}}).
\end{equation*}
\end{theorem}
\parindent=0mm\textbf{Proof}. Using the fact that the operators $L_{n_1,n_2; l_1, l_2}^{\alpha_1,\alpha_2;\beta_1,\beta_2}(f;q_{n_1}, q_{n_2};x, y)$ are linear
and positive together with the property $(ii)$ of the modulus of continuity, we obtain
\begin{eqnarray*}
\vert L_{n_1,n_2; l_1, l_2}^{\alpha_1,\alpha_2;\beta_1,\beta_2}(f;q_{n_1}, q_{n_2};x, y)-f(x, y)\vert \leq
\vert L_{n_1,n_2; l_1, l_2}^{\alpha_1,\alpha_2;\beta_1,\beta_2}(\vert f(t, s)-f(x,y)\vert;q_{n_1}, q_{n_2};x, y)-f(x, y)\\
\leq \omega (f;(\delta_{n_1}(x))^{\frac{1}{2}},(\delta_{n_2}(y))^{\frac{1}{2}}))\left(L_{n_1, l_1}^{\alpha_1, \beta_1}(1;q_{n_1},x)
+\frac{1}{(\delta_{n_1}(x))^{\frac{1}{2}}} L_{n_1, l_1}^{\alpha_1, \beta_1}(\vert t-x\vert;q_{n_1},x)\right) \\
\times\left(L_{n_2, l_2}^{\alpha_2, \beta_2}(1;q_{n_2},y)
+\frac{1}{{(\delta_{n_2}(y))^{\frac{1}{2}}}} L_{n_2, l_2}^{\alpha_2, \beta_2}(\vert s-y\vert;q_{n_2},y)\right).
\end{eqnarray*}
On applying the Cauchy-Schwartz inequality,
\begin{eqnarray*}
\vert L_{n_1,n_2; l_1, l_2}^{\alpha_1,\alpha_2;\beta_1,\beta_2}(f;q_{n_1}, q_{n_2};x, y)-f(x, y) \vert
\leq  \omega (f;(\delta_{n_1}(x))^{\frac{1}{2}},(\delta_{n_2}(y))^{\frac{1}{2}})) \\
\times\left(1+\frac{1}{(\delta_{n_1}(x))^{\frac{1}{2}}} (L_{n_1, l_1}^{\alpha_1, \beta_1}((t-x)^2;q_{n_1},x))^{\frac{1}{2}}\right) \\
\times\left(1+\frac{1}{{(\delta_{n_2}(y))^{\frac{1}{2}}}} (L_{n_2, l_2}^{\alpha_2, \beta_2}((s-y)^2;q_{n_2},y))^{\frac{1}{2}}\right),
\end{eqnarray*}
we obtain the required result.

\section{Approximation results}
In this section we prove some theorems regarding the degree of approximation for the bivariate operators through the Lipschitz class.
The Lipschitz class for the bivariate case, denoted by $Lip_M(\alpha_1, \alpha_2)$ is defined as under: \\
 Let $0< \alpha_1, \alpha_2\leq1$. A function $f$ is said to be in the class $Lip_M(\alpha_1, \alpha_2)$ if it satisfies the following inequality:
\begin{equation*}
\vert f(x, y)-f(x', y')\vert\leq M\vert x-x'\vert^{\alpha_1} \vert y-y'\vert^{\alpha_2}
\end{equation*}
for all $(x, y),(x', y')\in I_1\times I_2$.
Now we have the following theorem:
\begin{theorem}
Let $f\in Lip_M(\alpha_1, \alpha_2 $. Then
\begin{equation*}
\vert L_{n_1,n_2; l_1, l_2}^{\alpha_1,\alpha_2;\beta_1,\beta_2}(f;q_{n_1}, q_{n_2};x, y)-f(x, y)\vert \leq M
\sqrt{(\delta_{n_1}(y))^{\alpha_1}}\sqrt{(\delta_{n_2}(y))^{\alpha_2}}
\end{equation*}
holds for all $(x,y) \in J^2$.
\end{theorem}
\parindent=0mm\textbf{Proof}.
Using the hypothesis, we can write
\begin{eqnarray*}
  \vert L_{n_1,n_2; l_1, l_2}^{\alpha_1,\alpha_2;\beta_1,\beta_2}(f;q_{n_1}, q_{n_2};x, y)-f(x, y)\vert
  &\leq&L_{n_1,n_2; l_1, l_2}^{\alpha_1,\alpha_2;\beta_1,\beta_2}(\vert f(t, s)-f(x, y)\vert;q_{n_1}, q_{n_2};x, y)  \\
   &\leq& M L_{n_1,n_2; l_1, l_2}^{\alpha_1,\alpha_2;\beta_1,\beta_2}(\vert t-x\vert^{\alpha_1}\vert s-y\vert^{\alpha_2};q_{n_1}, q_{n_2};x, y)\\
 &=& M L_{n_1; l_1}^{\alpha_1;\beta_1}(\vert t-x\vert^{\alpha_1};q_{n_1};x) L_{n_2; l_2}^{\alpha_2;\beta_2}(\vert s-y\vert^{\alpha_2};q_{n_2};y).
\end{eqnarray*}
Now we apply the $H\ddot{o}lder's$ inequality with $u_1= \frac{2}{\alpha_1}, v_1= \frac{2}{2-\alpha_1}$ and $u_2 =\frac{2}{\alpha_2} , v_2
=\frac{2}{2-\alpha_2}$
to get
\begin{eqnarray*}
 \vert L_{n_1,n_2; l_1, l_2}^{\alpha_1,\alpha_2;\beta_1,\beta_2}(f;q_{n_1}, q_{n_2};x, y)-f(x, y)\vert
 &\leq& M L_{n_1; l_1}^{\alpha_1;\beta_1}((t-x)^2;q_{n_1};x)^{\frac{\alpha_1}{2}}
 L_{n_1; l_1}^{\alpha_1;\beta_1}(1;q_{n_1};x)^{\frac{2-\alpha_1}{2}} \\
 &{}&L_{n_2; l_2}^{\alpha_2;\beta_2}((s-y)^2;q_{n_2};y)^{\frac{\alpha_2}{2}}
 L_{n_2; l_2}^{\alpha_2;\beta_2}(1;q_{n_2};y)^{\frac{2-\alpha_2}{2}} \\
 &=& M (\delta_{n_1}(x))^{\frac{\alpha_1}{2}}(\delta_{n_2}(y))^{\frac{\alpha_2}{2}}\\
&=& M \sqrt{(\delta_{n_1}(y))^{\alpha_1}}\sqrt{(\delta_{n_2}(y))^{\alpha_2}}
\end{eqnarray*}
and this completes the proof of the theorem. \\

In the ensuing, we will use the following notations:
\begin{equation*}
C^1(I_1\times I_2) = \left\{f\in C(I_1\times I_2): f_x', f_y' \in C(I_1\times I_2)\right\}.
\end{equation*}
Now we will prove the following theorem:
\begin{theorem}
Let $f\in C^1(I_1\times I_2)$ and $(x, y) \in J^2$. Then
\begin{equation*}
\vert L_{n_1,n_2; l_1, l_2}^{\alpha_1,\alpha_2;\beta_1,\beta_2}(f;q_{n_1}, q_{n_2};x, y)-f(x, y)\vert \leq \Vert f_x'\Vert_{C(I_1\times I_2)}
\sqrt{\delta_{n_1}(x)}+ \Vert f_y'\Vert_{C(I_1\times I_2)}
\sqrt{\delta_{n_2}(y)}.
\end{equation*}
\end{theorem}
\parindent=0mm\textbf{Proof}.
For a fixed $(x, y) \in J^2$, let us write
\begin{equation*}
f(t, s)-f(x, y)= \int_{x}^{t}f_u'(u, s)d_qu +\int_{y}^{s}f_v'(x, v)d_qv.
\end{equation*}
Operating $L_{n_1,n_2; l_1, l_2}^{\alpha_1,\alpha_2;\beta_1,\beta_2}$ on both sides of the above equation, we get
\begin{eqnarray*}
\vert L_{n_1,n_2; l_1, l_2}^{\alpha_1,\alpha_2;\beta_1,\beta_2}(f;q_{n_1}, q_{n_2};x, y)-f(x, y)\vert&\leq&
L_{n_1,n_2; l_1, l_2}^{\alpha_1,\alpha_2;\beta_1,\beta_2}\left(\vert \int_{t}^{x}\vert f_u'(u, s)\vert d_qu\vert;q_{n_1}, q_{n_2};x, y \right)\\
&+& L_{n_1,n_2; l_1, l_2}^{\alpha_1,\alpha_2;\beta_1,\beta_2}\left(\vert \int_{y}^{s}\vert f_v'(x, v)\vert d_qv\vert;q_{n_1}, q_{n_2};x, y \right).
\end{eqnarray*}
Now since
\begin{equation*}
\vert \int_{t}^{x}\vert f_u'(u, s)\vert d_qu\vert \leq \Vert f_x'\Vert_{C(I_1\times I_2)} \vert t-x\vert,
\end{equation*}
and
\begin{equation*}
\vert \int_{y}^{s}\vert f_v'(x, v)\vert d_qv\vert \leq \Vert f_y'\Vert_{C(I_1\times I_2)} \vert s-y\vert,
\end{equation*}
we obtain
\begin{eqnarray*}
  \vert L_{n_1,n_2; l_1, l_2}^{\alpha_1,\alpha_2;\beta_1,\beta_2}(f;q_{n_1}, q_{n_2};x, y)-f(x, y)\vert &\leq&\Vert f_x'\Vert_{C(I_1\times I_2)}
  L_{n_1, l_1}^{\alpha_1,\beta_1}\left(\vert t-x\vert; q_{n_1};x \right)
  +\Vert f_y'\Vert_{C(I_1\times I_2)}L_{n_2, l_2}^{\alpha_2,\beta_2}\left(\vert s-y\vert; q_{n_2};y \right).\\
\end{eqnarray*}
Using the Cauchy-Schwartz inequality, we obtain
\begin{eqnarray*}
\vert L_{n_1,n_2; l_1, l_2}^{\alpha_1,\alpha_2;\beta_1,\beta_2}(f;q_{n_1}, q_{n_2};x, y)-f(x, y)\vert&\leq& \Vert f_x'\Vert_{C(I_1\times I_2)}
 \left(L_{n_1, l_1}^{\alpha_1,\beta_1}\left((t-x)^2; q_{n_1};x \right)\right)^{\frac{1}{2}}
 \left(L_{n_1, l_1}^{\alpha_1,\beta_1}\left(1; q_{n_1};x \right)\right)^{\frac{1}{2}} \\
 &+&\Vert f_y'\Vert_{C(I_1\times I_2)}
 \left(L_{n_1, l_1}^{\alpha_1,\beta_1}\left((s-y)^2; q_{n_2};y \right)\right)^{\frac{1}{2}}
 \left(L_{n_1, l_1}^{\alpha_1,\beta_1}\left(1; q_{n_2};y \right)\right)^{\frac{1}{2}} \\
  &=& \Vert f_x'\Vert_{C(I_1\times I_2)}\sqrt{\delta_{n_1}(x)}+\Vert f_y'\Vert_{C(I_1\times I_2)}\sqrt{\delta_{n_2}(y)}.
\end{eqnarray*}
Hence the theorem. \\

Now we define the following:\\
If $f\in C(I_1\times I_2), ~\delta>0$, then the partial moduli of continuity of $f$ with respect to $s$ and $t$, is defined by
\begin{equation*}
\tilde{\omega_1}(f; \delta) = sup\big\{\vert f(x_1, t)-f(x_2, t)\vert: t\in I_2 ~and~ \vert x_1-x_2\vert \leq \delta\big\}
\end{equation*}
and
\begin{equation*}
\tilde{\omega_2}(f; \delta) = sup\big\{\vert f(s, y_1)-f(s, y_2)\vert: s\in I_1 ~and~ \vert y_1-y_2\vert \leq \delta\big\}.
\end{equation*}
Overtly they satisfy the properties of the usual modulus of continuity. \\

Next we have the following theorem:
\begin{theorem}
Let $f\in C(I_1\times I_2)$ and $(x, y)\in J^2$. Then
\begin{equation*}
\vert L_{n_1,n_2; l_1, l_2}^{\alpha_1,\alpha_2;\beta_1,\beta_2}(f;q_{n_1},q_{n_2};x, y)-f(x, y)\vert\leq2\tilde{\omega_1}(f;\sqrt{\delta_{n_1}(x)})+
2\tilde{\omega_1}(f;\sqrt{\delta_{n_1}(y)})
\end{equation*}
holds.
\end{theorem}
\parindent=0mm\textbf{Proof}.
Making use of the definition of partial moduli of continuity, we obtain
\begin{eqnarray*}
  \vert L_{n_1,n_2; l_1, l_2}^{\alpha_1,\alpha_2;\beta_1,\beta_2}(f;q_{n_1},q_{n_2};x, y)-f(x, y)\vert &\leq &
   L_{n_1,n_2; l_1, l_2}^{\alpha_1,\alpha_2;\beta_1,\beta_2}(\vert f(t,s)-f(x,y)\vert;q_{n_1},q_{n_2};x, y) \\
  &\leq&  L_{n_1,n_2; l_1, l_2}^{\alpha_1,\alpha_2;\beta_1,\beta_2}(\vert f(t,s)-f(t,y)\vert;q_{n_1},q_{n_2};x, y) \\
  &+&L_{n_1,n_2; l_1, l_2}^{\alpha_1,\alpha_2;\beta_1,\beta_2}(\vert f(t,y)-f(x,y)\vert;q_{n_1},q_{n_2};x, y) \\
  &\leq& \tilde{\omega_1}(f;\delta_{n_1}(x))\big(L_{n_1,l_1, }^{\alpha_1,\beta_1}(1;q_{n_1};x)+ \frac{1}{\sqrt{\delta_{n_1}(x)}}
  L_{n_1,l_1, }^{\alpha_1,\beta_1}(\vert t-x\vert;q_{n_1};x)\big)\\
  &+&\tilde{\omega_2}(f;\delta_{n_2}(y))\big(L_{n_2,l_2}^{\alpha_2,\beta_2}(1;q_{n_2};y)+ \frac{1}{\sqrt{\delta_{n_2}(y)}}
  L_{n_2,l_2}^{\alpha_2,\beta_2}(\vert s-y\vert;q_{n_2};y)\big)
\end{eqnarray*}
On using the Cauchy-Schwartz inequality, we obtain
\begin{eqnarray*}
\vert L_{n_1,n_2; l_1, l_2}^{\alpha_1,\alpha_2;\beta_1,\beta_2}(f;q_{n_1},q_{n_2};x, y)-f(x, y)\vert &\leq & \tilde{\omega_1}(f;\delta_{n_1}(x))
\big(1+ \frac{1}{\sqrt{\delta_{n_1}(x)}}L_{n_1,l_1, }^{\alpha_1,\beta_1}((t-x)^2;q_{n_1};x)\big) \\
&+& \tilde{\omega_2}(f;\delta_{n_2}(y))
\big(1+ \frac{1}{\sqrt{\delta_{n_2}(y)}}L_{n_2,l_2 }^{\alpha_2,\beta_2}((s-y)^2;q_{n_2};y)\big) \\
&=& 2\tilde{\omega_1}(f;\sqrt{\delta_{n_1}(x)})+2\tilde{\omega_2}(f;\sqrt{\delta_{n_1}(y)}),
\end{eqnarray*}
and the proof is completed. \\

\end{document}